\documentclass[12pt]{amsart} \usepackage{amsmath,amsthm,amsfonts,amssymb, amscd, wasysym}

\usepackage{amssymb,eucal,mathrsfs,mathbbol}
\usepackage[all]{xy} 
\usepackage{amssymb, amsthm,
amsmath} \usepackage[english]{babel} \usepackage{amsfonts} \usepackage{tikz}

\newtheorem{theorem}{Theorem}[section]

\theoremstyle{definition} 
 
 \newtheorem{prop}[theorem]{Proposition}
\newtheorem{defn}[theorem]{Definition} \newtheorem{cor}[theorem]{Corollary}

 \newtheorem{rem}[theorem]{Remark}
\newtheorem{lem}[theorem]{Lemma} \theoremstyle{remark}

 \newcommand{\ra}{\rightarrow} 
 
\newcommand{\Z}{\mathbb Z} \newcommand{\ot}{\otimes} \newcommand{\mtc}{\mathcal}
 \newcommand{\Lam}{\Lambda} 
 \newcommand{\al}{\alpha} \newcommand{\eps}{\epsilon}
\newcommand{\bn}{\begin}  
 \newcommand{\D}{\Delta} 
  
\newcommand{\rh}{\rightharpoonup} \newcommand{\lh}{\leftharpoonup}
\numberwithin{equation}{section} 

\newcommand{\dw}{\downarrow}
\newcommand{\uw}{\uparrow}  \newcommand{\bl}{\begin{lem}}
\newcommand{\nc}{\newcommand} \nc{\el}{\end{lem}} 
\newcommand{\ovr}{\overline}
 \newcommand{\ch}{\chi} 
\newcommand{\mtr}{\mathrm} 
\newcommand{\beq}{\begin{equation}} \newcommand{\eeq}{\end{equation}}
\newcommand{\ncm}{\newcommand} 
\ncm{\hsa}{Hopf subalgebra of } \ncm{\ses}{semisimple} \ncm{\x}{$} \ncm{\bwt}{\bowtie}
\ncm{\mi}{\mtr{I}} \ncm{\cZ}{\mtc{Z}}\ncm{\xra}{\xrightarrow}
\ncm{\cb}{\mtc{B}}\ncm{\ca}{\mtc{A}} \ncm{\irr}{\Irr} \ncm{\Irr}{\mathrm{Irr}}
\ncm{\co}{\mtc{O}}
\ncm{\blue}{\textcolor[rgb]{.00, .00, 1.00}} \ncm{\red}{\textcolor[rgb]{1.00, .00, .00}}
\ncm{\md}{\medbreak} \ncm{\green}{\textcolor[rgb]{.00, 1.00, .00}}
\ncm{\Gm}{\Gamma}\ncm{\ind}{\mtr{Ind}}\ncm{\res}{\mtr{Res}}
\numberwithin{equation}{section} \ncm{\bq}{\beq}\ncm{\mto}{\mapsto}\ncm{\opl}{\oplus}
\ncm{\eq}{\eeq}
\ncm{\sg}{\sigma}
\nc{\bpf}{\begin{proof}} \nc{\epf}{\end{proof}}\nc{\bt}{\begin{theorem}}
\nc{\et}{\end{theorem}}
\newcommand{\dz}{\diez}
\ncm{\bp}{\begin{prop}} \ncm{\ep}{\end{prop}}

\ncm{\np}{\newpage} \ncm{\ebl}{\end{thebibliography}} \ncm{\bbl}{\begin{thebibliography}}
\ncm{\chd}{_{ _{\ch}}} \ncm{\ald}{_{ _{\al}}} \ncm{\cP}{\mathcal{P}} \ncm{\ei}{e_i}
\ncm{\eij}{e_{i,\;j}} \ncm{\bne}{\begin{enumerate}}
\ncm{\ene}{\end{enumerate}}\ncm{\bdef}{\begin{defn}} \ncm{\edf}{\end{defn}}
\ncm{\stab}{\mtr{Stab}} \ncm{\bc}{\begin{cor}}

\ncm{\ec}{\end{cor}} \ncm{\er}{\end{rem}} \ncm{\br}{\begin{rem}}

\ncm{\beqn}{\begin{equation*}} \ncm{\eeqn}{\end{equation*}}

\ncm{\bd}{\begin{document}} \ncm{\ed}{\end{document}}

\ncm{\cm}{\mathcal{M}} \ncm{\rep}{\mtr{Rep}} \ncm{\btw}{\bowtie} \ncm{\cd}{\mtc{D}}
\ncm{\cop}{\mtr{cop}}

\ncm{\bea}{\begin{eqnarray}} \ncm{\eea}{\end{eqnarray}} \ncm{\beanon}{\begin{eqnarray*}}
\ncm{\eeanon}{\end{eqnarray*}}\ncm{\ek}{\eps|_K}\ncm{\diez}{\#}

\ncm{\cC}{\mtc{C}}

\ncm{\cc}{\mtc{C}}

\ncm{\HKer}{\mtr{HKer}} \ncm{\LKER}{\mtr{LKER}} \ncm{\aad}{\mtr{ad}} \ncm{\Dr}{\mtr{D}}
\ncm{\cD}{\mathcal{D}} \ncm{\G}{\mathcal{G}} \ncm{\Dc}{\mtc{D}} \ncm{\E}{\mtc{E}}
\ncm{\fp}{\mtr{FPdim}} \ncm{\Vc}{\mtr{Vec}} \ncm{\cK}{\mtc{K}} \ncm{\cM}{\mtc{M}}
\ncm{\cE}{\mtc{E}} \ncm{\cS}{\mtc{S}} \ncm{\End}{\mtr{End}} 
\ncm{\ce}{\mtc{E}}\ncm{\inv}{\mtr{Inv}}\ncm{\blu}{\blue}\ncm{\lb}{\label}\ncm{\nfc}{normal
fusion subcategory} \ncm{\bxt}{\boxtimes} \ncm{\cl}{\mtc{L}} \ncm{\cbk}{\mtc{K}}
\ncm{\cbx}{\mtc{X}} 
\title[Semisimple Hopf
algebras] {Normal Hopf subalgebras of semisimple Drinfeld doubles} \author{Sebastian 
Burciu} \address{Inst.\ of Math.\ ``Simion Stoilow" of the Romanian Academy P.O. Box
1-764, RO-014700, Bucharest, Romania}

\email{sebastian.burciu@imar.ro, smburciu@gmail.com} 
\begin{document} \subjclass[2000]{Primary
16W30, 18D10} \keywords{Drinfeld doubles; fusion categories; normal Hopf subalgebras; abelian extensions} 
\begin{abstract} A description of all normal Hopf subalgebras of a semisimple Drinfeld double is given. This is obtained by considering  an analogue of Goursat's lemma concerning fusion subcategories of Deligne products of two fusion categories. As an application we show that the Drinfeld double of any abelian extension is also an abelian extension.
\end{abstract} 
\maketitle
\section*{Introduction and the main results}
Let $A$ be a semisimple Hopf algebra and $D(A)$ be its Drinfeld double. It is well  known that the category $\mtr{Rep}(D(A))$ of finite dimensional representations of $D(A)$
is a modular tensor category which is braided equivalent to the Drinfeld center of the fusion category $\mtr{Rep}(A)$ (see \cite{BaKi} or \cite{kas}). From this
point of view Drinfeld doubles form a special class of quasitriangular Hopf algebras that
play a very important role in the classification of semisimple Hopf algebras. Another
important notion in the classification of semisimple Hopf algebras is that of a normal
Hopf subalgebra of a Hopf algebra. For example, the class of semisolvable Hopf algebras
introduced in the paper \cite{MW} is constructed starting from a  tower of normal Hopf subalgebras.
\md In this paper we will give a general description for all normal Hopf subalgebras of
$D(A)$ in terms of normal Hopf subalgebras of $A$ and $A^*$ and some intermediate group data. In particular we show that if $A$ is simple as Hopf algebra then all normal Hopf subalgebras of $D(A)$ are central group subalgebras. \md  In order to show that the Drinfeld double of an abelian extension is a group theoretical Hopf algebra in \cite[Theorem 1.3]{nat-kac} it is shown that this Drinfeld double is equivalent to an $R$-twist of the (twisted) Drinfeld double of a finite group. We show in Theorem \ref{mac} that the Drinfeld double of any abelian extension is also an abelian extension. As another application we obtain a
description of all minimal normal Hopf subalgebras of a semisimple Drinfeld double $D(A)$. 
\md In order to obtain the description mentioned above we make use of the theory of fusion
categories, in particular  we apply a highly nontrivial quantum analogue of Goursat's
Lemma for groups. This quantum analog result first appeared in the unpublished manuscript \cite{dg} and it is presented
here in Theorem \ref{dg}. \md Using this result, a categorical description for all Hopf
subalgebras of $D(A)$ is presented. This enables us to construct a new class of Hopf
subalgebras of $D(A)$ in Theorem \ref{hs}. This class is parameterized by two normal Hopf
subalgebras $L_1, L_2$ of $A$ and a finite group $G$ with some additional properties. The
corresponding Hopf subalgebra of $D(A)$ is denoted by $B(L_1, L_2, G)$. It is shown in Theorem \ref{gennormal} that
any Hopf subalgebra of $D(A)$ with normal intersections with $A$ and $A^*$ is of the
type $B(L_1, L_2, G)$ mentioned above. \md  Necessary and sufficient conditions for $B(L_1, L_2, 1)$ to be a normal Hopf
subalgebra of $D(A)$ are given in Theorem \ref{comm}. 
Moreover, with these conditions satisfied it is shown that the quotient Hopf algebra $D(A)//B(L_{1}, L_{2}, 1)$ is a bicrossed product of the Hopf algebras $L_{1}^{*\cop}$ and $A//L_{2}$.  
Along the way we also obtain some
other new results concerning Hopf subalgebras of a semisimple Drinfeld double. For example
in Theorem \ref{trivh} it is shown that any Hopf subalgebra of $D(A)$ having trivial
intersections with $A$ and $A^*$ is a group algebra. \md This paper is organized as
follows. The first section recalls few basic results on fusion categories and semisimple
Hopf algebras that are needed through the rest of the paper. Section \ref{deligne}
presents the quantum analogue of Goursat's lemma that appeared first in \cite{dg}.  In the
next section the Hopf subalgebras $B(L_1, L_2, G)$ of $D(A)$ are constructed. Section
\ref{nm} shows that any Hopf subalgebra of $D(A)$ with normal intersections with $A$ and
$A^*$ is of the form $B(L_1, L_2, G)$. In the last section some examples and applications
of the above results are presented. For any abelian extension $A$ the Appendix describes the group structures from Theorem \ref{mac} that realizes $D(A)$ as an abelian extension.\md We work over an algebraic closed field $k$ of
characteristic zero. We use a short version  $\D(x)=x_1\ot x_2$ of Sweedler's notation for comultiplication. All the other Hopf algebra notations of this paper are similar to those
used in \cite{Montg}.

\section{Normal fusion subcategories}\lb{nrm}
In this
section we recall few basic facts on fusion categories and normal fusion subcategories
from \cite{GN} and \cite{dg} that are needed through the rest of the paper.
\subsection{General conventions on fusion categories} As usually, by a fusion category we
mean a $k$-linear semisimple rigid tensor category $\cc$ with finitely many isomorphism
classes of simple objects, finite dimensional spaces of morphisms, and such that the unit
object of $\cc$ is simple. We refer the reader to \cite{ENO} for a general theory of such
categories. \md Let $\mtc{C}$ be a fusion category and denote by $\co(\mtc{C})$ its set of
simple objects considered up to isomorphism.  Recall that the Grothendieck ring $K_0(\cc)$
of $\cc$ is the free $\Z$-module generated by the isomorphism classes of simple objects of
$\cc$ with the multiplication induced by the tensor product in $\cc$. The Grothendieck
ring $K_0(\cc)$ is a based unital ring  (see for example \cite{GN} for definition of based
rings). The isomorphism classes $[X]$ of simple objects of $\cc$ form a $\Z$-basis for
$K_0(\cc)$. \md Let $L_{[X]}$ be the linear operator given by left multiplication by $[X]$
in the based ring $K_0(\cc)$. Then the Frobenius-Perron dimensions of an object $X\in\cc$
is defined as the largest positive eigenvalue (Frobenius-Perron eigenvalue) of the matrix
associated to $L_{[X]}$ with respect to the linear basis given by $\co(\cc)$ of
$K_0(\cc)$. This eigenvalue is usually denoted by $\fp(X)$. \md A fusion subcategory $\cd$
of a fusion category $\cc$ is a full abelian replete subcategory $\cd$ of $\cc$ which is
closed under the tensor product. \subsection{Pointed fusion categories} A fusion category
$\cc$ is called pointed if any simple object of $\cc$ is invertible. It is well known that
pointed fusion categories are of the form $\mtr{Vec}^{\omega}_{G}$ for a finite group $G$
and some three cocycle $\omega \in H^3(G, k^*)$. Moreover the category of representations
$\rep(H)$ of a semisimple  Hopf algebra is pointed if and only $H$ is commutative, i.e
$H=kG^*$ as Hopf algebra. In this case one has $ \rep(H)=\mtr{Vec}_{G}$, i.e. $\omega=1$.
\subsection{Universal gradings offusion categories} 
Let $\mtc{C}$
be a fusion category. A grading of a fusion category $\mtc{C}$ by a group $G$ (or a
$G$-grading) is a map $\mtr{deg} : \mtc{O}(\mtc{C}) \ra G$ with the following property:
for any simple objects $X, Y, Z \in \mtc{ C}$ such that $X \otimes Y $ contains $Z$ one
has $\mtr{deg}(Z) = \mtr{deg}(X) \mtr{deg}(Y)$ . This  corresponds to a decomposition
$\mtc{C} = \oplus_{ g\in G} \mtc{C}_g$, where $\mtc{C}_g\subset \mtc{C}$ is the full
additive subcategory generated by the simple objects of $\cc$ with degree $g$. The abelian
subcategory $\mtc{C}_1$ is a fusion subcategory of $\cc$, called the trivial component of
the grading.\md Recall that a $G$-grading is said to be trivial if $\mtc{C}_1 = \mtc{C}$.
It is also said to be faithful if the map $\mtr{deg} : \mtc{O}(\mtc{C}) \ra G$ is
surjective. For any fusion category $\mtc{C}$, as explained in \cite [Sect. 3.2]{GN},
there is a notion of universal grading  whose grading group is called the universal
grading group of $\cc$ and it is denoted by $U_{ _{\mtc{C}}}$. Then the trivial component
of the universal grading is the fusion subcategory $\mtc{C}_{ _{ad}}$. Recall that $\mtc{
C}_{ _{ad}}$ is defined as the smallest fusion subcategory of $\mtc{C}$ containing all
$X\ot X^*$ for each simple object $X \in \co(\cc)$ of $\mtc{C}$. The universality of this
grading consists in the fact that any other grading factors through the universal grading.
\md The following Lemma appears in \cite{dg} and it will be needed in the proof of Theorem
\ref{dg}.

\bn{lemma} \label{subcad} Let $\mtc{C}$ be a fusion category and $$ \mtc{C}=\oplus_{g\in
U_{ _{\mtc{C}}}} \mtc{C}_g $$ be its universal grading. There is a one-to-one
correspondence between fusion subcategories $\mtc{D} \subset \mtc{C}$ containing
$\mtc{C}_{ _{ad}}$ and subgroups $G \subset U_{\mtc{C}}$, namely $$ D \mapsto  G_{
_{\mtc{D}}} := \{g \in U_{\mtc{C}}| \mtc{D} \cap \mtc{C}_g \neq 0 \} $$ and $$ G \mapsto
\mtc{D}_G :=\oplus _{g\in G}\; \mtc{C}_g.$$ \end{lemma}
\subsection{Few more results on normal fusion categories} \subsection{The commutator subcategory of a
fusion category} Recall the notion of commutator subcategory from \cite{GN}. If $\mtc{D}$
is a fusion subcategory of $\mtc{C}$ then $\mtc{D}^{co}$ is the full abelian subcategory
of $\mtc{C}$ generated by those objects $X$ such that $X\otimes X^* \in \mtc{D}$. If
$K_0(\cc)$ is a commutative ring then it is shown in \cite{GN} that $\cd^{co}$ is a fusion
subcategory of $\cc$. In general the result is not true. \subsection{Definition of a
normal fusion subcategory } The following notion of a normal tensor functor was introduced
in \cite{brn}. If $\cc$ and $\ce$ are fusion categories then a tensor functor $F:\cc\ra
\ce$ is normal if the following property is satisfied: if $m_{ \cc}(1_{\cc},\;F(X))>0$ for
some simple object $X \in \irr(\cc)$ then $F(X)=\fp(X)1_{\cc}$. Recall the multiplicity
form $m_{ \cc}$ defined on $K_0(\cc)$ by $m_{ \cc}([X], [Y])=\delta_{ [X], [Y]}$ for any
two simple objects $X, Y \in \co(\cc)$.

The fusion subcategory of $\cc$ generated by all objects $X\in \cc$ with 
$F(X)=\fp(X)1_{\cc}$ is called the kernel of $F$ and denoted by $\ker_{\cc} F$. Recall from 
\cite{brn} that a fusion subcategory $\cd \subset \cc$ is called normal if there is a
normal tensor functor $F:\cc\ra \ce$ such that $\cd=\ker_{\cc} F$. \ncm{\rad}{\mtr{rad}}
\br\lb{int} It is easy to check that if $\cd$ is a normal fusion subcategory of $\cc$ then
$\cd\cap \cc_1$ is a normal fusion subcategory of $\cc_1$ for any fusion subcategory
$\cc_1$ of $\cc$. \er \subsection{On the commutator and the radical of a fusion
subcategory} In this subsection we recall few notions and results from \cite{bb} that will
be used in this paper. For a fusion subcategory $\cd$ of a fusion category $\cc$ define
its radical as \bq \mtr{rad}_{\cc}(\cd)=\{X \in \co(\cc) \;|\; X^{\ot n} \in \cd \text {
for some } n \geq 1\} \eq If $K_0(\cc)$ is commutative then clearly $\rad_{\cc}(\cd)$ is a
fusion subcategory.\md Note that for any fusion category $\cd$ one has $\cd^{co}\subseteq
\rad_{\cc}(\cd)$ since $X^*$ always appear as a constituent of a tensor power of $X$ (see
also \cite[Remark 3.2]{nat}). In \cite{bb} it is shown that for a normal fusion
subcategory $\cd$ the above radical and the commutator coincide. \md For a fusion category
$\cc$ denote by $\inv(\cc)$, the set of all invertible objects (up to isomorphism) of $\cc$.
Suppose that $\cd$ is a normal fusion subcategory of $\cc$ fitting into an exact sequence
\beq \cd\ra\cc \xra{F} \ce \eeq of fusion categories. Then by \cite[Proposition 5.3]{bb}
it follows that $\cd^{co}$ is a fusion subcategory of $\cc$ and one has that \beq\lb{comb}
\cd^{co}=\mtr{rad}_{\cc}(\cd)=\{X\in \cc\;\;|\;\; F(X)=\fp(X)M, M \in \inv(\ce)\;\}. \eeq

\bn{example}  \lb{commh} Let $L$ be a normal Hopf subalgebra of $A$. Then the commutator
ideal $[A,L]$ is a Hopf ideal of $A$ and it is shown in \cite[Theorem 2.2]{cejm2} that
\begin{equation}\label{co} \mtr{Rep}(A//L)^{co}=\mtr{Rep}(A/[A,L]). \end{equation}
\end{example} 

The following Corollary is well known. For the sake of completeness we
include a sketch proof below. \bc\lb{vec} Any group extension of $\mtr{Vec}$ is a pointed
fusion category. \ec \bpf Let $ \cc=\oplus_{g \in G}\cc_g $ with the trivial component
$\cc_1=\mtr{Vec}$. By \cite[Proposition 5.4]{bb} it follows that $\cc_g$ has up to
isomorphism just one simple object $X_g$ which is also invertible. \epf \br\lb{radgr} Note
that if $ \cc=\oplus_{g \in G}\cc_g $ is a graded fusion category then
$\cc_1^{co}=\rad_{\cc}(\cc_1)=\cc$. Indeed any simple object $X \in \cc_g$ satisfies that
$X^{\ot\;n}\in \cc_{g^n}=\cc_1$ if $n$ is the order of $G$. This shows that
$\rad_{\cc}(\cc_1)=\cc$. On the other hand since $X^* \in \cc_{g^{-1}}$ it follows that
$X\ot X^* \in \cc_1$ and therefore $\cc_1^{co}=\cc$. \er \br\lb{costs} The previous
Proposition also shows that $\cd^{co}=\cc$.  Note that in the case of the universal
grading of $\cc=\rep(A)$ one has the above situation since $\cc_1=\rep(A//K(A))$ is a
normal fusion subcategory of $\rep(A)$.\er \subsection{Normal Hopf subalgebras and normal
fusion subcategories of $\mtr{Rep}(A)$} Let $A$ be a finite dimensional semisimple Hopf
algebra over $\mathbb{C}$. Then $A$ is also cosemisimple \cite{Lard}. The character ring
$C(A)$ of $A$ is a semisimple subalgebra of $A^*$ \cite{Z} and it has a vector space basis
given by the set $\mtr{Irr}(A)$ of irreducible characters of $A$. Moreover,
$C(A)=\mtr{Cocom}(A^*)$, the space of cocommutative elements of $A^*$. By duality, the
character ring of $A^*$ is a semisimple subalgebra of $A$ and $C(A^*)=\mtr{Cocom}(A)$. If
$M$ is an $A$-representation with character $\chi$ then $M^*$ is also an
$A$-representation with character $\chi^*=\chi \circ S$. This induces an involution
$``\;^*\;":C(A)\ra C(A)$ on $C(A)$. Let $m_{ _A}(\ch,\;\mu)$ be the multiplicity form on
$C(A)$. We use the notation $G(A)$ for the set of grouplike elements of $A$. 
\md
It is well known that $\mtr{Rep}(A)$ is a fusion category. Moreover there is a
maximal central Hopf subalgebra $K(A)$ of $A$ such that $\mtr{Rep}(A)_{
_{ad}}=\mtr{Rep}(A//K(A))$, see \cite{GN}. Since $K(A)$ is commutative it follows that
$K(A)=k[U_{ _A}]^*$ where $U_{ _A}$ is the universal grading group of
$\mtr{Rep}(A)$.

\bn{example}If $A=kG$ then $K(A)=k\mtr{Z}(G)$ and $U_{
_A}=\widehat{\mtr{Z}(G)}$, the linear dual group of the center $\mtr{Z}(G)$ of
$G$.\end{example} \md Let $\mtc{D}$ be a fusion subcategory of $\mtr{Rep}(A)$ and
$\mtc{O}(\mtc{D})$ be its set of simple objects. Then $I_{ _{\mtc{D}}}:=\cap_{V \in\mtc{O}(\mtc{D})}\mtr{Ann}_A(V)$ is a Hopf ideal in $A$ \cite{PQ} and
$\mtc{D}=\mtr{Rep}(A/I_{ _{\mtc{D}}})$. For a fusion category $\mtc{D} \subset
\mtr{Rep}(A)$ define its regular character as $r_{ _{\mtc{D}}}:=\sum_{X \in
\co(\mtc{D})}\dim_{k}(X)\ch_X$ where $\mtr{Irr}(\mtc{D})$ is the set of irreducible objects of $\mtc{D}$ and $\ch_X$ is the character of $X$ as $A$-representation. Thus $r_{
_{\mtc{D}}}\in C(A)$. \br\label{main} Let $A$ be a finite dimensional semisimple Hopf
algebra and $\mtc{D}$ be a fusion subcategory of $\mtr{Rep}(A)$. Then  \cite[Lemma
1]{Masnr} implies that $\mtc{D}=\mtr{Rep}(A//L)$, for some normal Hopf subalgebra $L$ of
$A$ if and only if the regular character $r_{ _{\mtc{D}}}$ of $\mtc{D}$ is central in
$A^*$. \er

\bn{defn} Let $L$ be a normal Hopf subalgebra of $A$. An irreducible character $\al$ of
$L$ is called $A$-stable if there is a character $\ch \in \mtr{Rep}(A)$ such that 
$\ch\dw^{ ^A}_{ _L}=\frac{\ch(1)}{\al(1)}\al$. Such a character $\ch \in \Irr(A)$ is said
to seat over the character $\al \in \Irr(L)$. \end{defn} The set of all irreducible
$A$-characters seating over $\al$ is denoted by $\Irr(A|_{\al})$. Denote by
$G_{A}^{st}(L)$ the set of all $A$-stable linear characters of $L$. Clearly $G_{A}^{st}(L)$ is a subgroup of the group of grouplike elements $G(L^*)$ of the dual
Hopf algebra of $L$. \md Proposition \ref{comm} and Example \ref{commh} imply the
following: \bl\label{power} Let $A$ be a semisimple Hopf algebra and $L$ be a normal Hopf
subalgebra of $A$.  Suppose that $M$ is an irreducible $A$-module affording the character
$\ch$. Then $M \in \mtc{O}(\mtr{Rep}(A//L)^{co})$ if and only if $L$ acts trivially on
some tensor power $M^{\ot\;n}$ of $M$. In these conditions $\ch\dw^{ ^A}_{ _L}=\ch(1)\psi$
for some $A$-stable linear character $\psi$ of $L$. \el \bpf Since $\rep(A//L)$ is a
normal fusion subcategory of $\rep(A)$ one has that
$(\mtr{Rep}(A//L)^{co})=\rad(\mtr{Rep}(A//L))$. This shows that $M \in
\mtc{O}(\mtr{Rep}(A//L)^{co})$ if and only if $L$ acts trivially on some tensor power
$M^{\ot\;n}$ of $M$. If $M\in \rad(\mtr{Rep}(A//L))$ then \cite[Proposition 5.3]{bb}
implies that $\ch\dw^{ ^A}_{ _L}=\ch(1)\psi$ for some $A$-stable linear character $\psi$
of $L$. If $M \in \rad(\mtr{Rep}(A//L))$ then by Equation \eqref{comb} for the functor
$F=\mtr{res}^A_L$ it follows that  $\ch\dw^{ ^A}_{ _L}=\ch(1)\psi$ for some $A$-stable
linear character $\psi$ of $L$. \epf For any character $\ch$ of $A$ denote by
$\mtr{Irr}(\ch)$ the set of all irreducible constituents of $\ch$. For $\cc\subset
\rep(A)$ we denote by $\cc\ch$ the full abelian subcategory of $\rep(A)$ generated by all
the objects $a\ch$ with $a \in \co(\cc)$. \bp\label{dual-costs} Suppose $A$ is a
semisimple Hopf algebra and $L$ is a normal Hopf subalgebra. Moreover suppose that there
is a  fusion subcategory $\cc \subset \rep(A)$ with a faithful $G$-grading $\cc=\oplus_{g
\in G}\cc_g$ whose trivial component is $\mtr{Rep}(A//L)$.  Then for any $\ch \in \cc_g$
one has that $\ch\dw^{ ^A}_{ _L}=\ch(1)\psi_g$ for some $A$-stable linear character
$\psi_g$ of $L$. Moreover the set $\co(\cc_g)=\Irr(A|_{\psi_g})$ and it can be described
as \beqn \cc_g=\cc_1\ch=\mtr{Irr}(\psi_g\uw_{ _L}^{ ^A}). \eeqn \ep \ncm{\onh}{On the
other hand\;} \bpf By Remark \ref{radgr} it follows that $\ch\in \rep(A//L)^{co}$.
Therefore by Lemma \ref{power} one has $\ch\dw^{ ^A}_{ _L}=\ch(1)\psi_g$ for some
$A$-stable linear character $\psi_g$ of $L$. \onh $\cc_1=\rep(A//L)$ is a normal fusion
subcategory of $\cc$ and it follows by \cite[Proposition 5.4]{bb} $\cc_g=\cc_1\ch$. The
equality $\cc_1\ch=\mtr{Irr}(\psi_g\uw_{ _L}^{ ^A})$ follows then by Frobenius
reciprocity. \epf \section{Fusion subcategories of a Deligne direct products of
categories}\lb{deligne} Let $\mtc{C}_1$ and $\mtc{C}_2$ be two fusion categories and form
the Deligne product $\cc:=\mtc{C}^1\boxtimes\mtc{C}^2$ . Identify $\cc_1$ and $\cc_2$ to
$\mtc{C}_1\boxtimes 1$ and $1 \boxtimes \mtc{C}_2$ respectively as fusion subcategories of
$\cc$. Recall that every simple object of $\cc$ is of the type $X_1\boxtimes X_2$ where
$X_i$ is a simple object of $\mtc{C}_i$. \md Let $\mtc{D} \subset \cc$ be any fusion
subcategory. Define $\mtc{L}^i(\mtc{D}):=\mtc{D}\cap \mtc{C}_i$, with $i=1,2$. Let also
$\mtc{K}_1(\mtc{D})$ be the fusion subcategory generated by all simple objects $X_1$ of
$\mtc{C}_1$ such that $X_1 \boxtimes X_2 \in \mtc{D}$ for some simple object $X_2$ of
$\mtc{C}_2$. Similarly define the fusion subcategory $\mtc{K}^2(\mtc{D})$. Clearly $
\mtc{L}^1(\mtc{D})\boxtimes\mtc{L}^2(\mtc{D})\subset
\mtc{D}\subset\mtc{K}^1(\mtc{D})\boxtimes\mtc{K}^1(\mtc{D})$. The following Theorem is an
analogue of Goursat's Lemma for direct product of groups and it appears in  \cite{dg}. For
the sake of completeness we include its proof below.

\bt\label{dg} Let $\mtc{D} \subset \mtc{C}_1\boxtimes\mtc{C}_2$ be a fusion subcategory.
Then there is a finite group $G$ and faithful $G$-gradings
$\mtc{K}_i(\mtc{D})=\bigoplus_{x \in G}\mtc{K}^i(\mtc{D})_x$ with trivial components
$\mtc{L}^i(\mtc{D})$, $i=1, 2$ such that \begin{equation}\label{81} \mtc{D}=\bigoplus_{x
\in G}\mtc{K}^1(\mtc{D})_x\boxtimes\mtc{K}^2(\mtc{D})_x. \end{equation} \et \ncm{\ac}{\{}
\ncm{\pd}{[} \bn{proof} First we will show that \begin{equation}\label{82} \mtc{D} \supset
\mtc{K}^1(\mtc{D})_{ _{ad}} \boxtimes \mtc{K}^2(\mtc{D})_{ _{ad}} = (\mtc{K}^1(\mtc{D}) 
\boxtimes \mtc{K}^2(\mtc{D}))_{ _{ad}} , \end{equation} where as above
$\mtc{K}^i(\mtc{D})_{ _{ad}}$ denotes the adjoint subcategory of $\mtc{K}^i(\mtc{D})$.
Note that if $X_1  \boxtimes  X_2 \in \mtc{O}(\mtc{D})$ then $(X_1 \otimes X^*_1 ) 
\boxtimes (X_2 \otimes X^*_2) \in \in \mtc{O}(\mtc{D})$. Since $X_i\otimes X^*_i$ contains
the unit object of $\cc_i$ it follows that $(X_1 \otimes X_1^*)\boxtimes 1 \in \mtc{D}$
and $1\boxtimes (X_2\otimes X^*_2) \in \mtc{D}$. Therefore $ \mtc{K}^i(\mtc{D})_{ _{ad}}
\subset \mtc{L}^i(\mtc{D})$\; \text{for} \;i=1,2. Since $\mtc{L}^1(\mtc{D})  \boxtimes
\mtc{L}^2(\mtc{D}) \subset \mtc{D}$, it follows that above inclusion implies Equation
\eqref{82}. Now let $U_{ _{\mtc{K}^i(\mtc{D})}}$ be the universal grading group of
$\mtc{K}_i(\mtc{D})$. Lemma \ref{subcad} applied to the inclusion Equation \eqref{82}
gives that \beqn \mtc{D} =\bigoplus_{\gamma \in \Gamma}(\mtc{K}^1(\mtc{D}) \boxtimes
\mtc{K}^2(\mtc{D}))_{\gamma} \eeqn for some subgroup $\Gamma \in U_{
_{\mtc{K}_1(\mtc{D})\boxtimes \mtc{K}^2(\mtc{D})}}= U_{ _{\mtc{K}^1(\mtc{D})}} \times U _{
_{\mtc{K}_2(\mtc{D})}}$. On the other hand it follows from the definition of
$\mtc{K}^i(\mtc{D})$ that the maps $\Gamma \ra U_{ _{\mtc{K}^1(\mtc{D})}}$ and $\Gamma \ra
U_{ _{\mtc{K}^2(\mtc{D})}}$ are surjective. Then Goursat's Lemma for groups (see
\cite{Gours}) implies that $\Gamma$ equals the fiber product $U_{
_{\mtc{K}^1(\mtc{D})}}\times_{ _{G}} U_{ _{\mtc{K}^2(\mtc{D})}}$ for some group $G$
equipped with group epimorphisms $U_{ _{\mtc{K}^i(\mtc{D})}} \ra G$, $i = 1, 2$. Then
these epimorphisms define faithful $G$-gradings of $\mtc{K}^i(\mtc{D})$ such that Equation
\eqref{81} holds. In this case note that $\mtc{L}^i(\mtc{D}) := \mtc{D} \cap \mtc{C}^i$
equals the trivial component of the $G$-grading of $\mtc{K}^i(\mtc{D})$ for all $1\leq i
\leq 2$. \end{proof} \br\lb{drad} By Remark  \eqref{radgr} it follows that in this
situation \beqn \rad(\mtc{L}^i(\cd)) \supseteq \mtc{L}^i(\cd)^{co} \supseteq
\mtc{K}_i(\cd) \eeqn for all $1 \leq i \leq 2$. \er \bp\lb{co} Let $\cd$ be any fusion
subcategory of $\cc_1\boxtimes \cc_2$. Then \bq \cd^{co}= \cl_1(\cd)^{co}\bxt
\cl_2(\cd)^{co} \eq and \bq \rad({\cd})=\rad(\cl_1(\cd)) \bxt \rad(\cl_2(\cd)) \eq \ep
\bpf Suppose that $X \bxt Y \in \cd^{co}$. Then $(X\ot X^*)\bxt (Y\ot Y^*) \in \cd$ and
therefore $X\ot X^* \in \cl_1(\cd)$ since $1_{\cc_2}$ is a constituent of $Y\ot Y^*$.
Similarly $Y\ot Y^*\in \cl_2(\cd)$. Thus $\cd^{co}\subset \cl_1(\cd)^{co}\bxt
\cl_2(\cd)^{co}$. On the other hand, clearly $\cl_1(\cd)^{co}\bxt \cl_2(\cd)^{co}\subset
\cd^{co}$. \md Suppose now that $X \bxt Y \in \rad(\cd)$. Then $X^n \bxt Y^n \in \cd$ for
some $n \geq 1$. As in \cite[Remark 3.2]{nat} there is $m \geq 1$ such that
$Y^{mn}=(Y^n)^m$ contains the unit element of $\cc_2$. This implies that $X^{mn}\in
\cl^1(\cd)$ and therefore $X \in \rad(\cl^1(\cd))$. Similarly $Y \in \rad(\cl^2(\cd))$
which shows that $\rad({\cd})\supset\rad(\cl_1(\cd)) \bxt \rad(\cl_2(\cd))$. The other
inclusion is immediate. \epf \bt\lb{trivf} Any fusion subcategory of $\cc_1\boxtimes
\cc_2$ with trivial intersections with $\cc_1$ and $\cc_2$ is pointed. \et \bpf With the
above notations it follows that $\mtc{L}_1(\cd)=\mtr{Vec}$ and $\mtc{L}_2(\cd)=\mtr{Vec}$.
Thus $\mtc{K}_i(\cd)$ are group extensions of $\mtr{Vec}$ and by Corollary \ref{vec} they
are pointed fusion subcategories. It follows that $\cd$ is also a pointed fusion category.
\epf \subsection{Normal fusion subcategories of a Deligne tensor product} Suppose now that
$\cd$ is a normal fusion subcategory of $\cc: = \cc_1 \bxt \cc_2$. Applying Theorem
\ref{dg} one obtains \beq \mtc{D}=\oplus_{x \in
G}\mtc{K}^1(\mtc{D})_x\boxtimes\mtc{K}^2(\mtc{D})_x \eeq where $\cbk^i(\cd)$ are faithful
graded fusion subcategories \beq \cbk^i(\cd)=\oplus_{x \in \cbx}\cbk^i(\cd)_x. \eeq
\ncm{\ck}{\mtc{K}} Moreover by Remark  \eqref{int}  one has that $\cl^i(\cd)$ are normal
fusion subcategories of $\cc^i$ and $\rad(\mtc{L}^i(\cd))=
\mtc{L}^i(\cd)^{co}\supseteq\mtc{K}_i(\cd)$ for all $1 \leq i \leq 2$.

\br Since $\cl^i(\cd)$ is a normal fusion subcategory of $\ck^i(\cd)$ it
follows by \cite[Proposition 5.4]{bb} that $\ck^i(\cd)_x=\cl^i(\cd)a^i_x$ for any simple
object with $a^i_x \in \ck^i(\cd)_x$. \er \section{Hopf subalgebras of semisimple Drinfeld
doubles}\lb{hsa} Recall that the Drinfeld double $D(A)$ of a Hopf algebra $A$ is defined
by $D(A) \cong A^{*cop} \otimes A$ as coalgebras with the multiplication given by \beq (g
\bowtie h)(f \bowtie l)=\sum g(h_1\rightharpoonup f \leftharpoonup S^{-1}h_3)\bowtie
h_2l,\eeq for all $h,l \in A$ and $f,g \in A^*$. Moreover its antipode is given by $S(f
\bowtie h)=S^{-1}(h)S(f)$. It is well known that $D(A)$ is a semisimple Hopf algebra if
and only if $A$ is a semisimple Hopf algebra \cite{Montg}. It is also known that $D(A)$ is
a cocycle twist of $A^{*\mtr{cop}}\otimes A$ and therefore
$\mtr{Rep}(D(A)^*)=\mtr{Rep}(A)^{\mtr{rev}}\boxtimes\mtr{Rep}(A^*)$ where $\mtr{Rep}(A)^{\mtr{rev}}:=\rep(A^{\mtr{op}})$.
\subsection{On Hopf subalgebras of Drinfeld doubles $D(A)$}\label{hsalg}

If $H$ is a Hopf subalgebra of $D(A)$ then clearly one has the following inclusion \beqn
\mtr{Rep}(H^*)\subset \mtr{Rep}(D(A)^*)= \mtr{Rep}(A)^{\mtr{rev}}\boxtimes \mtr{Rep}(A^*).
\eeqn Let $\cc_1:=\mtr{Rep}(A)^{\mtr{rev}}$, $\cc_2:=\mtr{Rep}(A^*)$ and $\cd:=\rep(H^*)$.
Theorem \ref{dg} implies that $\mtr{Rep}(H^*)=\bigoplus_{x \in G}\mtc{K}^1(\cd)_x
\boxtimes \mtc{K}^2(\cd)_x $. Moreover one has the inclusions $\mtc{K}^1(\cd)_{
_{ad}}\subset \mtc{L}^1(\cd)\subset \mtc{K}^1(\cd) \subset \mtr{Rep}(A)^{rev}, $ and
$\mtc{K}^2(\cd)_{ _{ad}}\subset \mtc{L}^2(\cd) \subset \mtc{K}^2(\cd)  \subset
\mtr{Rep}(A^*)$ \ncm{\cx}{\mtc{X}} with the faithful  $G$-gradings $
\mtc{K}^i(\cd)=\bigoplus_{x \in G}(\mtc{K}^i(\cd))_x $. Note that in this situation
$\mtc{K}^i(\cd)_1=\mtc{L}^i(\cd)$ for all $i=1,2$. Also one can write
$\mtc{L}^1(\cd)=\rep(M_1^*)^{\mtr{rev}}$ where $M_1:=H\cap A^*$  is a Hopf subalgebra of
$A^*$. Similarly $\mtc{L}^2(\cd)=\rep(M_2^*)$ where $M_2:=H\cap A$ is a Hopf subalgebra of
$A$. \bt\lb{trivh} Any Hopf subalgebra of $D(A)$ with trivial intersections with $A$ and
$A^*$ is a group algebra. \et 

\bpf Let $H$ be a such Hopf
subalgebra of $D(A)$. It follows that $\rep(H^*)$ regarded as above, as a fusion subcategory of
$\mtr{Rep}(A)^{\mtr{rev}}\boxtimes \mtr{Rep}(A^*) $, has trivial intersections with both
$\mtr{Rep}(A)^{\mtr{rev}}$ and $\mtr{Rep}(A^*)$. Then Theorem \eqref{trivf} implies that
$\rep(H^*)$ is pointed. Thus $H^*=kG^*$ and $H=kG$. \epf 

For any right $A$-comodule $M$
with comodule structure $\rho:M \ra A\ot M$ denote by $C_M$ the coalgebra of
coefficients of $M$. Recall that $C_{M}$ is the smallest subcoalgebra $C\subset A$ with the property
that $\rho(M)\subset C\ot M$ \cite{Lar}. If $d \in C(A^*)$ is the character of $M$ as an
$A$-comodule then $C_M$ is also denoted by $C_d$.\md By duality any left $A$-module $V$
with associated character $\ch \in C(A)$ can be regarded as a right $A^*$-module and one
can associate to it as above its subcoalgebra of coefficients $C_{\ch}\subset A^*$. \md
Recall from \cite{NZ} that a subset $X \subset \mtr{Irr}(H^*)$ is closed under
multiplication if for every two elements $c, d \in X$ in the decomposition of the product
$cd=\sum_{e \in \mtr{Irr}(H^*)}m_{c,d}^{e}e$ one has $e\in X$ whenever $m_e\neq 0$. Also a
subset $X \subset \mtr{Irr}(H^*)$ is closed under $``\;^*\;"$ if $x^* \in X$ for all $x
\in X$.

\br\lb{closed} Following \cite [Theorem 6]{NZ} it follows that any subset $X \subset
\Irr(H^*)$ closed under multiplication and $``\;^*\;"$ generates a Hopf subalgebra $H(X)$ of
$H$ defined by \beq H(X):=\oplus_{x \in X}C_x. \eeq \er \subsection{On a class of Hopf
subalgebras of $D(A)$} Suppose that $L_1$ and $L_2$ are two normal Hopf subalgebras of $A$
and let $G$ be a finite group that can be simultaneously embedded in  $G^{st}_A(L_1)$ and
$ G^{st}_{A^*}((A//L_2)^*)$ via the emebddings $\psi_1:G \hookrightarrow G^{st}_A(L_1)$
and respectively $\psi_2:G \hookrightarrow G^{st}_{A^*}((A//L_2)^*)$. Let $B(L_1, L_2,G)$
be the subcoalgebra of $D(A)$ defined by \beq B(L_1, L_2,G)=\bigoplus_{x \in
G}C_{\psi_1(x)\uw^{ ^{A}}_{ _{L_1}}}\bowtie C_{ \psi_2(x)\uw^{ ^{A^*}}_{ _{(A//L_2)^*}}}.
\eeq 
\bp\lb{hs} Suppose that $L_1$ and $L_2$ are two normal Hopf subalgebras of
$A$ and let $G$ be a finite group as above. Then $B(L_1, L_2, G)$ is a Hopf subalgebra of
$D(A)$. \ep \bpf By Remark \eqref{closed} it is enough to show that the irreducible
comodule characters of $B(L_1,L_2,G)$ form a set closed under tensor product and duality.
Note that an irreducible comodule character of the Hopf subalgebra $B(L_1, L_2,G)$ is of
the type $\ch\bwt d$ with $\ch \in \Irr(A)$ seating over $\psi_1(x)$ (i.e. $\ch \in
\Irr(A|_{\psi_1(x)})$) and $d \in \Irr(A^*)$ seating over $\psi_2(x)$ (i.e. $d\in
\Irr(A^*|_{\psi_2(x)})$), for some $x \in G$. Suppose moreover that $\ch' \in \Irr(A)$
seats over $\psi_1(y)$ and $d' \in \Irr(A^*)$ seats over $\psi_2(y)$ for some other $y \in
G$. It follows that for the product in $\rep(D(A)^*)$ \beq (\ch\bwt d)(\ch'\bwt
d')=\ch\ch'\bwt dd' \eeq the irreducible constituents of $\ch\ch'$ are all seating over
$\psi_1(xy)$ and the irreducible constituents of $dd'$ are all seating over $\psi_2(xy)$.
This shows that the above set is closed under product. Moreover the set is closed under
duality since $\ch^*\in \Irr(A|_{\psi_1(x^{-1})})$ and $d^*\in
\Irr(A^*|_{\psi_2(x^{-1})})$. \epf

Note that one has \beq\lb{inta} B(L_1, L_2,G)\cap A=L_2 \; \text{and}\;  B(L_1, L_2,G)\cap
A^*=(A//L_1)^*.\eeq \md If $G=1$  then $\psi_1$ and $\psi_2$ are trivial and we denote the
Hopf subalgebra  $B(L_1, L_2, \{1\})$ simply by $B(L_1,L_2)$.
 \subsection{Hopf
subalgebras of $D(A)$ with normal intersections with $A$ and $A^*$}
\bt\label{gennormal} Let $A$ be a semisimple Hopf algebra. Then any Hopf subalgebra $H$ of
$D(A)$ with normal intersections with $A$ and $A^*$ is of the form $B(L_1, L_2, G)$
defined above. \et \bn{proof} Let $\cd:=\rep(H^*)\subset \rep(A)^{rev}\boxtimes
\rep(A^*)$. Using the notations from subsection \eqref{hsalg} it follows that
$\mtc{L}^1(\cd)$ and $\mtc{L}^2(\cd)$ are normal fusion subcategories of $\rep(A)^{rev}$ and
respectively $\rep(A^*)$.  Therefore using Remark \ref{main} it follows that $M_1$
is a normal Hopf subalgebra of $A^*$ and one may suppose that $M_1:=(A//L_1)^*$ for some 
normal Hopf subalgebra $L_1$ of $A$. Thus $\mtc{L}^1(\cd)=\rep(A//L_1)^{rev}$. Similarly
$\mtc{L}^2(\cd)=\rep(L_2^*)$ for some normal Hopf subalgebras $L_1, L_2$ of $A$. Suppose
that $\eta \bowtie d \in \mtc{K}^1(\cd)_x \boxtimes \mtc{K}^2(\cd)_x \subset \cd$ for a
given $x \in G$. Then Lemma \ref{dual-costs} implies that $\eta \dw_{ _{L_1}}^{
_A}=\eta(1)f_x$ for some $A$-stable linear character $f_x$ of $L_1$. This shows that $G$
can be regarded as a subgroup of $G^{st}_{ A}(L_1)$ via the map $\psi_1$ given by $x
\mapsto f_x$. By duality, the same argument applied for $\mtc{K}^2(\cd)_x$ gives that $d
\dw^{ ^{A^*}}_{ _{(A//L_2)^*}}=\eps(d)g_x$ for some $A^*$-stable linear character $g_x$ of
$(A//L_2)^*$. Therefore $G$ can also be identified with a subgroup of $G^{st}_{
A^*}((A//L_2)^*)$ via the map $\psi_2$ given by $x \mapsto g_x$. Moreover by the same
Proposition \ref{dual-costs} the set of simple objects of $\ck^1(\cd)_x$ can be identified
with the set of irreducible modules of $\psi_1(x)\uw^A_{L_1}$. A similar result holds for
$\ck^2(\cd)_x$ and the proof of the Theorem is complete.\end{proof} \bc\label{cocor} With
the notations from Theorem \ref{gennormal} it follows that \beqn \mtr{Rep}(H^*)\subset
\mtr{Rep}(A//L_1)^{co}\boxtimes \mtr{Rep}(L_2^*)^{co} \eeqn \ec \bn{proof} It follows
directly for above by applying Remark \ref{power}. \end{proof} \section{Normal Hopf
subalgebras of Drinfeld doubles $D(A)$}\lb{nm} In this section we give the description for
all normal Hopf subalgebras of Drinfeld doubles $D(A)$. \subsection{ Commutativity between
$L_1$ and $L_2$} We need the following lemma: \bl\lb{commute} Suppose that $K$ is a Hopf
subalgebra of $H$ and $x \in H$. Then $ Sm_2xm_1=\eps(m)x $ for all $m \in K$ if and only
if $mx=xm$ for all $m \in K$. \el \bpf One has that $xm=m_1S(m_2)xm_3=m_1\eps(m_2)x=mx$.
The converse also follows since $Sm_2xm_1=Sm_2m_1x=\eps(m)x$. \epf \br Suppose that $K$ is
a normal Hopf subalgebra of a Hopf algebra $H$ and let $A=H//K$ be the quotient Hopf
algebra of $H$ via $\pi_K:H\ra H//K$. Then since $\pi^* :A^* \ra H^*$ is an injective Hopf
algebra map it follows that $\pi^*(A^*)$ is a  Hopf subalgebra of $H^*$. Under this
identification it can be checked that \begin{equation}\label{forml} \pi^*(A^*)=\{f\in
H^*|f(hm)=f(h)\eps(m)\;\; ; h\in H,\;m\in K\} \end{equation} is the set of linear
functionals invariant under left (right) $K$-translations. \er For two Hopf subalgebras
$L_1, L_2$ of $A$ denote as usually by $[L_1, L_2]$ the commutator ideal generated by
$ab-ba$ with $a\in L_1$ and $b \in L_2$. Note that two Hopf subalgebras $L_1$ and $L_2$ of
$A$ commute elementwise if and only if $[L_1, L_2]=0$. \md Let $K$ be a Hopf subalgebra of
$A$ and $\Lam_K\in A$ be its idempotent integral. Then it is well known (see \cite{Bker})
that the induced representation $A\ot_Kk$ is isomorphic to $A\Lam_K$ via the map
$a\ot_K1\mapsto a\Lam_K$. \md Note that if $K$ is a normal Hopf subalgebra of $A$ then
$\Lam_K$ is a central element of $A$  by \cite[Lemma 1]{Masnr}. Then $K$ acts trivially on
the induced module $A\ot_Kk$. Indeed, since the induced module is isomorphic to $A\Lam_K$
one has that $ma\Lam_K=m\Lam_Ka=\eps(m)a\Lam_K$ for all $m \in K$ and $a \in A$. Corrolary 2.5 
of \cite{Bker} shows that $K$ is a normal Hopf subalgebra of $A$ if and only if $K$ acts
trivially on the induced module $A\ot_Kk$

\bn{theorem}\label{commutwise} Let $A$ be a semisimple Hopf algebra and $L_1, L_2$ be two
normal Hopf subalgebras of $A$. With the above notations if $B(L_1, L_2, G)$ is a normal
Hopf subalgebra of $D(A)$ then one has $[L_1, L_2]=0$ and $[(A//L_1)^*,\;(A//L_2)^*]=0$
\end{theorem}

\bn{proof} Let $K:=B(L_1, L_2, G)$. If $K$ is a normal Hopf subalgebra of $D(A)$ then as
above one has $m((f \bowtie b)\ot_K1)=\eps(m)((f \bowtie b)\ot_K1)$ for all $m \in K$ and
any $f \in A^*$, $b \in A$. In particular if $m \in L_2$ then one has \begin{eqnarray*} %
m((f \bowtie b)\ot_K1) &=&
f_1(Sm_3)f_3(m_1)f_2 \bowtie m_2b\ot_K1 \\ & = & f_1(Sm_3)f_3(m_1)f_2 \bowtie
b_1(Sb_2m_2b_3)\ot_K1\\ &=&f_1(Sm_3)f_3(m_1)f_2 \bowtie b_1\ot_K(Sb_2m_2b_3)1 \\ & = &
f_1(Sm_2)f_3(m_1)f_2 \bowtie b\ot_K1 \end{eqnarray*}

This implies that $f_1(Sm_2)f_3(m_1)f_2 \bowtie b\ot_K1 =\eps(m)f \bowtie b\ot_K1$ and the
previous remark gives that \begin{equation}\label{cond} (f_1(Sm_2)f_3(m_1)f_2 \bowtie
b)\Lam_K =\eps(m)(f \bowtie b)\Lam_K \end{equation} Note that \beqn \Lam_K=\sum_{x \in
G}\psi_1(x)\uw^A_{L_1}\bwt \psi_2(x)\uw^{A^*}_{(A//L_2)^*} \eeqn for all $f \in A^*$ and
$b \in A$. \md Since for $x \neq y$ one has that $C_{\psi_2(x)\uw^{A^*}_{(A//L_2)^*}}\cap
C_{\psi_2(y)\uw^{A^*}_{(A//L_2)^*}}=0$ there is a functional $g \in A^*$ such that
$g=\eps$ on the coalgebra $L_2=C_{\psi_2(1)\uw^{A^*}_{(A//L_2)^*}}$ and $g$ is zero on all
the other coalgebras $C_{\psi_2(x)\uw^{A^*}_{(A//L_2)^*}}$ with $x \neq 1$. \md Then put
$b=1$ in Equation \eqref{cond} and apply $\mtr{id}\otimes g$ on both therms. One obtains
that $f_1(Sm_2)f_3(m_1)f_2t_{A//L_1}=\eps(m)ft_{A//L_1}$ for all $f \in A^*$. Evaluating
both sides of the above equality at any $l \in L_1$ it follows that $Sm_2lm_1=\eps(l)m$
since $t_{(A//L_1)}(l)=\eps(l)$. Then Lemma \ref{commute} shows that $lm=ml$ for all $m
\in L_2$ and any $l \in L_1$. Thus $L_1$ and $L_2$ commute elementwise. \md It is known
that $D(A) \cong D(A^{*\;op\;cop})^{op}$ as Hopf algebras via $f \bwt a \mapsto a \bwt f$
(see for instance \cite[Theorem 3]{minimal}). Note that under the above isomorphism one has that $B(L_1,
L_2, G)=B((A//L_2)^*, (A//L_1)^*, G)$ as a Hopf subalgebra of $D(A^{*\;op\;cop})$. Indeed
note that $(A//L_1)^*\bwt L_2 =(A^*//(A//L_2)^*)^*\bwt (A//L_1)^*$ inside
$D(A^{*\;op\;cop})$. \md Then the above argument of this proof applied to $A^{*\;op\;cop}$
implies that $B((A//L_2)^*, (A//L_1)^*, G)$ is closed under the left adjoint action of
$A^*$ if and only if $[(A//L_2)^*,\;(A//L_1)^*]=0$. \end{proof} \subsection{Normal Hopf
subalgebras of type $B(K,L)$} One more preliminary lemma is needed. \bl \lb{adjd}Let $K$
be a normal Hopf subalgebra of a semisimple Hopf algebra $A$. Then
$m(g(a_2)Sa_3a_1)=(g(a_2)Sa_3a_1)m$ for all $m \in K$, $g \in (A//M)^*$ and any $a \in A$.
\el \bpf Since $(A//K)^*$ is a normal Hopf subalgebra of $A^*$ it follows that
$Sf_1gf_2\in (A//K)^*$ for all $f \in A^*$. Equation \eqref{forml} implies \beqn
<Sf_1gf_2, am>=\eps(m)<Sf_1gf_2, a> \eeqn for all $a \in A$ and all $m \in K$. It is easy
to see that the above equality can be written as
$f(Sm_1Sa_1g(a_2m_2)a_3m_3)=\eps(m)f(Sa_1g(a_2)a_3)$. Since $f$ was chosen arbitrary and
$g \in (A//K)^*$ it follows that \beqn Sm_1(Sa_1g(a_2)a_3)m_2=\eps(m)(Sa_1g(a_2)a_3) \eeqn
for all $a \in A$ and all $m \in K$. Applying previous lemma it follows that \beqn
m(Sa_1g(a_2)a_3)=m(Sa_1g(a_2)a_3) \eeqn for all $a \in A$ and any $m \in K$. Now the
result follows by passing to $K^{cop} \subset A^{\cop}$ as a normal Hopf subalgebra. \epf
For a trivial group, denote the Hopf subalgebra $B(K,L, 1)$ by $B(K,L)$. We now can prove
the following theorem:

 \bt\label{comm} If $K$ and $L$ are normal Hopf subalgebras of $A$
then $B(K,L)$ is a normal Hopf subalgebra of $D(A)$ if and only if $[K,\;L]=0$ and
$[(A//K)^*,\;(A//L)^*]=0$. \et

\bn{proof}  Note that $B(K,L):=(A//K)^{*\cop}\bowtie L$. If $B(K,L)$ is a normal Hopf
subalgebra of $D(A)$ then by Theorem \ref{commutwise} it follows that the two commutator
ideals vanish. \md Conversely, suppose that $[K,\;L]=0$ and $[(A//K)^*,\;(A//L)^*]=0$. \md
It will be shown first that if  $[K,\;L]=0$ then the Hopf subalgebra $B(K,L)$ is closed
under the adjoint action of $A^{*\cop}$ regarded as a Hopf subalgebra of $D(A)$. Note that
under the adjoint action of $A^{*\cop}$ one has that \beqn f.(g\bwt l)=f_2(g \bwt
l)Sf_1=f_2g(l_1\rh Sf_1\lh Sl_3)\bwt l_2, \eeqn for all $f \in A^{*\cop}$, $g \in
(A//K)^{*\cop}$ and $l \in L$.

Thus in order to show that $f.(g\bwt l) \in B(K, L)$ it is enough to show that $f_2g(a\rh
Sf_1\lh b)\in (A//K)^*$ for all $a, b \in L$.  \md By Equation \eqref{forml} $f_2g(a\rh
Sf_1\lh b) \in (A//K)^{*\cop}$  if and only  if $$<f_2g(a\rh Sf_1\lh b),
\;cm>=\eps(m)<f_2g(a\rh Sf_1\lh b),\;c>$$ for all   $f \in A^{*\cop}$, $g \in
(A//K)^{*\cop}$, $c \in A$ and $m \in K$. Thus one has to show that: \beqn<f,
SaSm_3Sc_3Sbc_1m_1g(c_2m_2)>=\eps(m)<f, SaSc_3Sbc_1g(c_2)>\eeqn for all $f \in A^*$, $g
\in (A//K)^{*\cop}$, $c \in A$, $a, b\in L$ and $m \in K$. Since $f \in A^*$ is arbitrary
this is equivalent to $SaSm_3Sc_3Sbc_1m_1g(c_2m_2)=\eps(m)SaSc_3Sbc_1g(c_2)$. Therefore it
is enough to show that \beq\lb{idc} Sm_3Sc_3Sbc_1m_1g(c_2m_2)=\eps(m)Sc_3Sbc_1g(c_2)\eeq
for all $g \in (A//K)^{*\cop}$, $c \in A$, $b\in L$ and $m \in K$. On the other hand since
$g\in (A//K)^{*\cop}$ the above equation can also be written as \beq\lb{idcc} Sm_2\pd
Sc_3Sbc_1g(c_2)]m_1=\eps(m)\pd Sc_3Sbc_1g(c_2)]\eeq for all $g \in (A//K)^{*\cop}$, $c \in
A$, $b\in L$ and $m \in K$. Thus applying Lemma \ref{commute} for $x=Sc_3Sbc_1g(c_2)$ the
above equation is equivalent to $m[Sc_3Sbc_1g(c_2)]=[Sc_3Sbc_1g(c_2)]m$  for all $a \in A$
and $x, y \in L$. \md On the other hand note that applying Lemma \eqref{adjd} one has
\begin{eqnarray*} m[Sc_3Sbc_1g(c_2)]&=&m[(Sc_5Sbc_4)(Sc_3g(c_2)c_1)]\\
&=&[(Sc_5Sbc_4)m(Sc_3g(c_2)c_1)]\\ &=& Sy(Sa_5xa_4)k(Sa_3g(a_2)a_1)\\ &=&
[Sc_3Sbc_1g(c_2)]m \end{eqnarray*} which shows that Equation \eqref{idcc} is satisfied.
\md One has $D(A) \cong D(A^{*\;op\;cop})^{op}$ as Hopf algebras via $f \bwt a \mapsto a
\bwt f$ by \cite[Theorem 3]{minimal}. Then as in the proof of the previous Theorem, it
follows by the same argument applied to $A^{*\;op\;cop}$ that $(A//K)^*\bwt L$ is closed
under the adjoint action of $A$ if and only if $[(A//K)^*,\;(A//L)^*]=0$. \end{proof}
\bc\lb{incl} With the above notations suppose that $B(L_1, L_2, G)$ is a normal Hopf
subalgebra of  $D(A)$. Then $B(L_1,\;L_2)$ is also a normal Hopf subalgebra of $D(A)$. \ec
\bpf It follows by Theorem \ref{comm} that the above conditions are satisfied. \epf
\md
\subsection{Bicrossed product of Hopf algebras}

Recall the notion of bicrossed product of two Hopf algebras. Let $H$ be a Hopf algebra and $A, L$ be two Hopf subalgebras of $H$. We say that $(A, L)$ is a factorization of $H$ if the multiplication
map $m: A\ot L \ra H$ is bijective. Note that in this case $A \cap  L = k$ and by \cite[Theorem 2.7.3]{majid} it follows that $H$ is a bicrossed product Hopf algebra of $A$ and $L$. It is also known that $D(A)$ is a bicrossed product of $A^{*\cop}$ and $A$, see e.g. \cite{majid}.
\bp \label{bicross}Let $A$ be a semisimple Hopf algebra and $B(L_1,\;L_2)$  a normal Hopf subalgebra of $D(A)$. Then $D(A)//B(L_{1}, L_{2})$ is a bicrossed product of $L_{1}^{*\cop}$ and $A//L_{2}$.
\ep
\bpf
Denote $D(L_{1}, L_{2}):=D(A)//B(L_{1}, L_{2})$. Since $B(L_{1}, L_{2})=(A//L_{1})^{*\cop}\bwt L_{2}$ one clealry has the following inclusions of Hopf ideals \\$D(A)B(L_{1}, L_{2})^{+}\supset AL_{2}^{+}$ and $D(A)B(L_{1}, L_{2})^{+}\supset\; A^{*\cop}(A//L_{1})^{*\cop\; +}$. These inclusion determine the following Hopf algebra emebddings:\\
$
A//L_{2}\hookrightarrow D(L_{1}, L_{2})$
 and 
$
 L_{1}^{*\cop}:=A^{*\cop}//(A//L_{1})^{*\cop}\hookrightarrow D(L_{1},L_{2}).
$
\md 
Moreover, the multiplication map 
$
 L_{1}^{*\cop}\bwt A//L_{2} \ra D(L_{1},L_{2})
$
 given by $\ovr{f}\ot  \ovr{a}\mapsto \ovr{(f \bwt a)}$ is clearly surjective. A dimension argument implies that this map is bijective and therefore $D(L_{1},L_{2})$ is a bicrossed product of $L_{1}^{*\cop}$ and $A//L_{2}$.\epf
\section{Examples and applications}\lb{exam} In this section we give some general examples
of normal Hopf subalgebras of $D(A)$ together with some applications. Recall that a Hopf
algebra $A$ is called simple if there are no proper normal Hopf subalgebras of $A$. \md
\bn{example}\label{ka} Let $A$ be a semisimple Hopf algebra and $K(A)$ its maximal central
Hopf subalgebra. Then $K(A)$ is a normal Hopf subalgebra of $D(A)$. Indeed in this case
one has that $K(A)=B(A, K(A))$ and the conditions of the above theorem are satisfied.
\end{example} \bn{example} Let $A$ be a semisimple Hopf algebra. Then for any central Hopf
subalgebra $L\subset K(A)$ it follows that $(A//L)^*\bowtie A$ is a normal Hopf subalgebra
of $D(A)$. Indeed in this case one has that $(A//L)^*\bowtie A=B(L, A)$ and the conditions
of the above theorem are satisfied. \end{example} 
\subsection{Drinfeld doubles of abelian extensions}
\bl \label{comd}Suppose that $K$ is a normal commutative Hopf subalgebra of a semisimple Hopf algebra $A$. Then $B(K,K)$ is a commutative Hopf subalgebra of $D(A)$.
\el
\bpf
One has to show that $mf=f\bwt m$  for all $m \in K$ and any $f \in (A//K)^{*}$. Note that using Equation \eqref{forml} one has
$mf=(Sm_3\rh f \lh m_{1})\bwt m_{2}=f\eps(m_{1})\eps(S(m_{3}))\bwt m_{2}=f\bwt m.$
\epf
Recall that a Hopf algebra $A$ is called an abelian extension if there is an exact sequence of Hopf algebras  \beq\label{see} k \ra k^G\ra A\xra{\pi} kF\ra k \eeq for some finite
groups $G$ and $F$.
Next it will be shown that the Drinfeld double of an abelian extension is also an abelian extension. Note that in \cite[Theorem 1.3]{nat-kac} in order to show that $\rep(D(A))$ is group-theoretical it is shown that the Drinfeld double of an abelian extension is an $R$-twist of a twisted Drinfeld double $D^{\omega}(\Sigma)$ for sum finite group $\Sigma$.
\bt\label{mac}
If $A$ is an abelian extension then $D(A)$ is also an abelian extension.
\et
\bpf
Suppose as above that $A$ is
obtained as an abelian extension: \beqn k \ra k^G\ra A\xra{\pi} kF\ra k \eeqn for some finite
groups $G$ and $F$. Then $B(k^G, k^G)$ is a normal Hopf subalgebra of $D(A)$. Indeed note
that $(A//k^G)^*\cong k^F$ is also a commutative Hopf subalgebra of $A^*$. 
One has the following extension
\bq
k \ra B(k^{G}, k^{G})\ra D(A) \xra{\delta} B\ra k
\eq 
where $B:=D(A)//B(k^{G}, k^{G})$ and $\delta$ is the canonical Hopf projection.
Note that $B(k^{G}, k^{G})=(k^{F})^{\cop}\bwt k^{G}$ as a Hopf subalgebra of $D(A)$. Moreover since $B(k^{G}, k^{G})$ is commutative by Lemma \ref{comd} it follows that $B(k^{G}, k^{G})=k^{F^{\mtr{op}}\times G}$ .

On the other hand by Proposition \ref{bicross} one has that $D(A)//B(K,K)\cong kG
\bwt kF$ is a bicrossed product Hopf algebras. Since both $kG$ and $kF$ are pointed it follows that their bicrossed product is also pointed and therefore a group algebra $kM$. 
\epf
The group $M$ from above is obtained as an exact factorization of the groups $G$ and $F$. Therefore $M=G\btw F$. In Theorem \ref{m} from Appendix we give a complete description of the group structure of $M$. As expected, it will be shown that $M$ is isomorphic to the group $\Sigma$ from \cite[Theorem 1.3]{nat-kac}.
\bt\lb{triv} Any normal Hopf
subalgebra of $D(A)$ with trivial intersections with $A$ and $A^*$ is a group algebra
which is central in $D(A)$. \et \bpf Let $H$ be a Hopf subalgebra of $D(A)$ with $H\cap
A=k$ and $H\cap A^*=k$. Then $L_1=A$ and $L_2=k$. Therefore using the notations from
Theorem \ref{gennormal} it follows that $G$ is a subgroup of $G(A^*)\cap G(A)$ and \beq
H=\bigoplus_{x \in G}kf_x\bwt kg_x. \eeq

Since $H$ is a normal Hopf subalgebra of $A$ it follows that $\Lam_H$ is a central element
of $D(A)$. Therefore $\Lam_H=\sum_{x \in G}f_x \bwt g_x$ is central in $D(A)$. Note that
the condition $f\Lam_H=\Lam_Hf$ implies that \beqn \sum_{x \in G}ff_x\bwt g_x=\sum_{x
\in G}f_xf(g_x^{-1}? g_x)\bwt g_x \eeqn  which shows that $f_x\bwt g_x$ commutes to any
element $f \in A^*$. Similarly $f_x\bwt g_x$ commutes to any element $a \in A$ and
therefore this is a central element of $D(A)$. \epf \bp\lb{containinga} Any normal Hopf
subalgebra of $D(A)$ containing $A$ is of the type $(A//L)^*\btw A$ with $L\subseteq
K(A)$. By duality any normal Hopf subalgebra of $D(A)$ containing $A^*$ is of the type
$A^*\btw K$ with $K$ a normal Hopf subalgebra of $A$ such that $(A/K)^*\subseteq K(A^*)$.
\ep \bpf If $A\subseteq H:=B(L_1, L_2, G)$ then $L_2=A$ and therefore $G=1$. Thus
$H=(A//L_1)^*\btw A$ and the condition $[L_1, A]=k$ implies that $L_1\subseteq K(A)$. \epf
\bc If $A$ is a simple Hopf algebra then any normal Hopf subalgebra of $D(A)$ is a central
group algebra. \ec \bpf Since $A$ is a simple Hopf algebra it follows that also $A^*$ is a
simple Hopf algebra. Thus any Hopf subalgebra of $D(A)$ has trivial intersections with $A$
and $A^*$ and therefore this is a central group algebra by Theorem \eqref{triv}. \epf 
\bp Any minimal normal Hopf subalgebra of $D(A)$ is of one of the following types:
\bne \item $B(K, L)$ for a  pair of normal Hopf subalgebras of $A$ satisfying $[K,L]=k$
and $[(A//K)^*,\;(A//L)^*]=k$. \item A central group subalgebra of $D(A)$. \ene
\ep \bpf Let $B(K, L,G) $ be a minimal normal Hopf subalgebra of $D(A)$. By Corollary
\ref{incl} it follows that $B(K,L)\subseteq B(K, L, G) $ is also a normal Hopf
subalgebra of $D(A)$. Therefore if $B(K,L)\neq B(K, L,G)$ then $B(K,L)=k$ which implies
that $K=A$ and $L=k$. It follows by Equation \eqref{inta} that $B(K, L, G) $ has trivial
intersections with $A$ and $A^*$. Thus by Theorem \ref{triv} it follows that $B(K, L,G)
$ is a central group subalgebra of $A^*$. \epf \bp If $A\cong K\otimes L$ as Hopf algebras
then $D(A)\cong B(K,L)\ot B(L,K)$ as Hopf algebras. \ep \bpf It is easy to see that
$A^*\cong (A//K)^*\otimes (A//L)^*$ as Hopf algebras. Then clearly $B(K,L)$ and $B(L,K)$
are normal Hopf subalgebras of $D(A)$. Morever $B(K,L)\cap B(L,K)=k$ and
$B(K,L)B(L,K)=D(A)$. It follows by \cite[Theorem 2.5]{cejm} that $D(A)\cong B(K,L)\ot
B(L,K)$ as Hopf algebras. \epf 
\section{Appendix: Drinfeld doubles of abelian extensions}
\subsection{Abelian extensions}
Let $\Sigma = FG$ be an an exact factorization of finite groups. This gives a right action  $\lhd : G \times F \ra G$ of $F$ on the set $G$, and a left action $\rhd : G \times F \ra F$ of $G$ on the set $F$ subject to the following two conditions:
\bq \label{c1}s\rhd xy = (s \rhd x)((s\lhd x) \rhd y)\;\;\;\;\;
 st\lhd x = (s\lhd (t \rhd x))(t \lhd x)\eq

The actions $\rhd$ and $\lhd$ are determined inside the group $\Sigma$ by the following relation 
\bq\label{m1}
gx = (g\rhd x)(g \lhd x)
\eq
for all $ x \in F$, $g \in G$. Note that $1 \rhd x=x$ and $s\lhd1=s$.
The quadruple $(G, F, \lhd, \rhd)$ is also called a matched pair of groups.

Consider the Hopf algebra $A=k^G \;^{\tau}\#_{\sg}\; kF$ from  \cite{Masext} which is a crossed product and coproduct formed using the above two actions. Recall that the structure of $A$ is given by:
\bq
(p_g\dz x)(p_h \dz y) = \delta_{g\lhd x,h} \sg_{g}(x, y)p_g \dz xy
\eq
\bq\label{del}
\D(p_g\dz x) = \sum_{st=g} p_s\#(t \rhd x) \otimes \tau_{x}(s,t)p_t\#x
\eq
where $p_{g}\in k^{G}$ is is the dual basis to the group element basis of $kG$. Without loss of generality we may suppose that $\tau_{x}(g, g^{-1})=1=\sg_{g}(x, x^{-1})$ for all $g\in G$ and $x \in F$ (see \cite[Lemma 3.6]{carp}). Moreover in this case the antipode of $A$ is given by
\bq
S(p_{g}\#x)=(g \rhd x)^{-1}p_{g^{-1}}=p_{(g\lhd x)^{-1}}\#(g \rhd x)^{-1}.
\eq 
Recall that here $\sg:F\ra F \times k^{G}$ is a normalized two cocyle written as $\sg(x,y):=\sum_{g\in G}\sg_{g}(x,y)p_{g}$
Similarly, $\tau:G\times G\ra k^{F}$ is a normalized two cocyle written as $\tau(g,h):=\tau_{x}(g, h)q_{x}$ where $q_{x} \in k^{F}$ is the dual basis to the group element basis of $kF$. 

Then the Hopf algebra $A$ from above fits into the abelian extension
\bn{equation}\begin{CD}\label{abext}k @ > >> k^G  @> i >>  A @ > \pi >> kF@ > >> k .\end{CD}\end{equation}

Note that under the identification $A=k^G \;^{\tau}\#_{\sg} kF$  one has that $\pi(p_{g }\#x)=\delta_{g,1}x$ and $i(p_{g})=p_{g}\#1$. This induces an isomorphism $A//k^{G}\xra{}kG$ via
$\overline{p_{a}\ot x}\mapsto \delta_{a,1}x$.
\md
Moreover, as it is shown in \cite{Masext} any abelian Hopf algebras fitting an exact sequence as in Equation \eqref{abext} is of the form presented above.
\subsection{Dual Hopf algebra of an abelian extension}The dual Hopf algebra $A^{*}$ fits the exact sequence 
\bn{equation}\begin{CD}\label{dabext}k @ > >> k^F  @> \pi^{*} >>  A @ > i^{*} >> kG@ > >> k \end{CD}\end{equation}
which shows that this is also an abelian extension.

Under the obvious identification $A^{*}\cong k^{F}
\;^{\sg^{*}}\#_{\tau^{*}} \;kG$ it can be shown that 
\bq
(q_{x}\ot a)(q_{y}\ot b)=\delta_{x, b \rhd y}\tau_{y}(a, b)q_{y}\dz ab
\eq
\bq\label{deld}
\D(q_{x}\ot a)=\sum_{uv=x}(q_{u}\ot a)\ot (\sg_{a}(u,v)q_{v}\ot (a \lhd u)).
\eq
Note that $A^{*}//k^{F}\xra{\cong}kG$ via $\ovr{q_{x}\ot a}\mapsto \delta_{x,1}a$
\subsection{On the Drinfeld double of abelian extensions} In this subsection we describe the group $M=G\bowtie F$ from the proof of Theorem \ref{mac}. More precisely we prove the following:
\bt\label{m}
Suppose that $A$ is an abelian extension fitting the exact sequence of Equation \eqref{abext}. Then $D(A)$ fits the following exact sequence
\bq
k \ra k^{F^{\mtr{op}}\times G}\ra D(A) \xra{\delta} k(G\bowtie F)\ra k
\eq
where the factorization of $G\bowtie F$ is given via the actions $
x \RHD a= (a^{-1} \lhd x^{-1})^{-1}$ and  $x \LHD a=(a^{-1}\rhd x^{-1})^{-1}$ for all $a \in G$ and $x \in F$. 
\et

For shortness of writing, for any $f \in A^{*}$ and any $a , b\in A$ in the proof below we denote the linear  functional $b \rh f \lh a \in A^{*}$ by $f(a?b)$. Note that $f(a?b)(x)=f(axb)$ for all $x \in A$.
\bpf Since $G^{op}\bowtie F$ is an exact factorization of groups it follows that $xa=(x \RHD a)(a \LHD x)$ for some matched pair of groups on $G^{\mtr{op}}$ and $F$.
\md
One has that $D(A)\cong A^{*}\bwt A$ as vector spaces and the canonical projection
\bq
D(A)\xra{\delta} k(G^{\mtr{op}}\bwt F)
\eq
from Theorem \ref{mac}
is given by $\ovr{q_{x}\ot b}\bwt \ovr{p_{a}\dz y}\mapsto \delta_{x,1}\delta_{a,1}b\bwt y$. Note that inside the quotient $k(G^{\mtr{op}}\bwt F)$ one has
$xa=\delta((\ovr{p_{1}\dz x})(\ovr{q_{1}\ot a}))$. On the other hand inside $D(A)$ it follows
\beanon
(p_{1}\dz x)(q_{1}\ot a)&=&<q_{1}\ot a, \; S(p_{1}\dz x)_{3}?(p_{1}\dz x)_{1}>(p_{1}\dz x)_{2}
\eeanon

Note that since we assumed $\tau_{x}(s, s^{-1})=1$ by Equation \eqref{del} one has 
\beanon
\D^{2}(p_{1}\dz x)&=&\sum_{s \in G}(p_{s}\dz (s^{-1}\rhd x))\ot \D(p_{s^{-1}}\dz x)\\&=& \sum_{m, n\in G}(p_{(mn)^{-1}}\dz (mn \rhd x)\ot (p_{m}\dz (n \rhd x))\ot \tau_{x}(m,n)(p_{n}\dz x).
\eeanon
This implies that
\beanon
(p_{1}\dz x)(q_{1}\ot a)&=&<q_{1}\ot a, \; S(p_{1}\dz x)_{3}?(p_{1}\dz x)_{1}>(p_{1}\dz x)_{2}
\\&=& \sum_{m, n\in G}<q_{1}\ot a, \;S(p_{n}\dz x)\tau_{x}(m,n)?(p_{(mn)^{-1}}\dz (mn \rhd x))>\bwt \\ & & \bwt (p_{m}\dz (n \rhd x)).
\eeanon
It follows that
\beanon
\delta((p_{1}\dz x)(q_{1}\ot a))&=&\sum_{m,n\in G}\delta(<q_{1}\ot a, \; S(p_{n}\dz x)?(p_{(mn)^{-1}}\dz (mn \rhd x))>)\times \\& \times &\delta((p_{m}\dz (n \rhd x)))\\&=& \sum_{n\in G}\delta(<q_{1}\ot a, \; S(p_{n}\dz x)?(p_{n^{-1}}\dz (n \rhd x))>)(n \rhd x)
\eeanon
Denote by $f^{n}_{x,a}:=<q_{1}\ot a, \; S(p_{n}\dz x)?(p_{(mn)^{-1}}\dz (mn \rhd x))>\;\in A^{*}$. Therefore the above formula can be written as
\bea\label{q1}
xa=\delta((p_{1}\dz x)(q_{1}\ot a))&=&\sum_{n \in G}\delta(f^{n}_{x,a})(n \rhd x)
\eea
Note that for any $f \in A^{*}$ one has \bq
f=\sum_{g \in G, \; x\in F}f(p_{b}\dz x)(q_{x}\ot b).
\eq  Therefore $\pi(f)=\sum_{b \in G}f(p_{b}\dz 1) b$. For $f^{n}_{x, a}$ it follows that
\beanon
\pi(f^{n}_{x, a})
&=&
\sum_{b \in G}f^{n}_{x,a}(p_{b}\dz 1) b
\\ &=&
 <q_{1}\ot a, \;(S(p_{n}\dz x))(p_{b}\dz 1)(p_{n^{-1}}\dz (n \rhd x))>b
\\ &=& 
 <q_{1}\ot a, \;(n \rhd x)^{-1}p_{n^{-1}}\dz (n \rhd x)>n^{-1}
\\&=& \delta_{b, n^{-1}}<q_{1}\ot a, \;p_{n^{-1}\lhd (n \rhd x)}\dz 1>n^{-1}
\\ &=& <q_{1}\ot a, \;p_{(n \lhd x)^{-1}}\dz 1>n^{-1}
\\ &=&
\delta_{a, (n \lhd x)^{-1}}n^{-1}
= \delta_{a^{-1}, (n \lhd x)}n^{-1}=\delta_{a^{-1} \lhd x^{-1}, n}(a^{-1} \lhd x^{-1})^{-1}
\eeanon
Thus Equation \eqref{q1} implies
\bq
xa= \sum_{n \in G}\delta(f^{n}_{x,a})(n \rhd x)= (a^{-1} \lhd x^{-1})^{-1} ((a^{-1} \lhd x^{-1})\rhd x)
\eq
This implies that
$
x \RHD a=( a^{-1} \lhd x^{-1})^{-1}
$
and 
$
x \LHD a=(a^{-1} \lhd x^{-1})\rhd x
$
It follows from Equation \eqref{c1} that $(a^{-1} \lhd x^{-1})\rhd x=(a^{-1}\rhd x^{-1})^{-1}$. Thus $x \LHD a=(a^{-1}\rhd x^{-1})^{-1}$.
\epf
\br
It follows that $M \cong \Sigma$ by the map $(a,x)\mapsto (x,a)$ since from Equation \ref{m1} the multiplication inside $\Sigma$ verifies
$xa=(a^{-1}x^{-1})^{-1}=(a^{-1}\rhd x^{-1})^{-1}(a^{-1}\lhd x^{-1})^{-1}$ for all $a \in G$ and all $x \in F$.
\er
\bibliographystyle{amsplain} 
\bibliography{nr-dr-db}
\end{document}